\documentclass[12pt,a4paper]{article}
\usepackage[T2A]{fontenc}
\usepackage[cp1251]{inputenc}
\usepackage{amsmath,amssymb,amsthm,amsfonts,amscd}
\begin{document}
\newtheorem{example}{Example}
\newtheorem{proposition}{Proposition}
\newtheorem{conjecture}{Conjecture}
\newtheorem{corollary}{Corollary}
\newtheorem{lemma}{Lemma}		
\newtheorem{remark}{Remark}
\newtheorem{theorem}{Theorem}
\newtheorem{definition}{Definition}
\newtheorem{problem}{Problem}
\def\K{{\bf K}}
\def\SH{{S\Bbb{H}}}
\def\SO{{S\Bbb{O}}}
\def\SU{{S\Bbb{U}}}
\def\S{{\Bbb{S}}}
\def\O{{\Bbb{O}}}
\def\R{{\Bbb R}}
\def\I{{\Bbb I}}

\def\Z{{\Bbb Z}}
\def\P{{\Bbb P}}
\def\RP{{\Bbb R}\!{\rm P}}
\def\N{{\Bbb N}}
\def\C{{\Bbb C}}
\def\H{{\Bbb H}}
\def\1{{\bf 1}}
\def\Q{{\bf Q}}
\def\A{{\bf A}}
\def\D{{\bf D}}
\def\k{{\bf k}}
\def\K{{\bf K}}
\def\B{{\bf B}}
\def\E{{\bf E}}
\def\F{{\bf F}}
\def\V{\vec{\bf V}}
\def\L{{\bf L}}
\def\G{{\bf G}}
\def\c{{\bf c}}
\def\i{{\bf i}}
\def\j{{\bf j}}
\def\f{{\bf f}}

\def\fr{{\operatorname{fr}}}
\def\st{{\operatorname{st}}}
\def\mod{{\operatorname{mod}\,}}
\def\cyl{{\operatorname{cyl}}}
\def\dist{{\operatorname{dist}}}
\def\grad{{\bf{grad}}}
\def\div{{\operatorname{div}}}
\def\rot{{\operatorname{rot}}}

\def\R{{\Bbb R}}
\def\B{{\bf B}}
\def\e{{\bf e}}
\def\L{{\bf L}}
\def\valpha{\vec{\alpha}}
\def\vxi{\vec{\xi}}
\sloppy
\title{Kairvaire Problems in Stable Homotopy Theory}
\author{Petr M. Akhmet'ev} 
\date{IZMIRAN RAS and MIEM HSE}
\maketitle
\medskip
\bigskip
\sloppy
\begin{abstract}

The Kervaire Problem is an unsolved problem in Stable Homotopy Theory.
The first interesting example is in dimension $30$. 
There exists a closed stably-parallelizable manifold
 $\tilde{M}^{30}$ with Arf-Kervaire invariant 1. It is unknown, if such a manifold exists in dimension 
$126$? The goal of the paper is to recall the construction of the manifold $M^{30}$ and investigate its properties.
\end{abstract}

\section*{Content}
\[  \]

1.  Kervaire Problem.

2.  Generalized Kervaire Problem.

3.  Strong Kervaire Problem.

4.  Barratt-Jones-Mahowald Theorem. James Theorem. Hill-Hopkins-Ravelel Theorem. Snaith Conjecture. 

5.  Positive solution of Strong Kervaire Problem in dimension $k=15$. 

6.  The Jones Manifold. 

7.  Browder Theorem. Eccles Theorem. The positive solution of Kervaire Problem in dimension $30$.

8.  The Eccles-Wood filtration. An estimation of the term of the filtration for the standard Kervaire $30$-manifold.

\section{Kervaire Problem}

Let  $f: M^{n-1} \looparrowright \R^n$ be an immersion of a closed, generally speaking, non-orientable manifold
in the codimension 1, assumed 
$n \equiv 4 \pmod{2}$. Let us consider the self-intersection manifold of $f$, denote this manifold by
$N^{n-2}$. Consider the characteristic class $w_2 \in H^2(N^{n-2};\Z/2)$, denote the characteristic number 
$ <w_2^{\frac{n-2}{2}};[N^{n-2}]> \in \Z/2$ by $\Theta(f)$. The most interesting case is $n=2^l-2$, $l \in \N, l \ge 2$.  

\begin{problem}{Kervaire Problem}\label{problem1}

Assume that $n=2^l-2$, $l \ge 2$. For which  $l$ the following equation $\Theta(f)=1$ for an immersion $f: M^{n-1} \looparrowright \R^n$ is possible? 
\end{problem}  

\section{Generalized Kervaire Problem}

Let $g: M^{n-k} \looparrowright \R^{n}$ be an immersion of a closed, generally speaking,non-oriented manifold of 
dimension $n-k$, where $n \equiv 2 \pmod{4}$, $n > k$. Denote by $\kappa: M^{k} \looparrowright \RP^{\infty}$ a prescribed characteristic class. In the case $n-k \equiv 1 \pmod{2}$ the characteristic class $\kappa$ is the oriented class over $M^{n-k}$: $\kappa = w_1(M^{n-k})$. The line bundle $\kappa^{\ast}(\gamma)$ ($\gamma$ is the only non-trivial bundle over $\RP^{\infty}$) is re-denoted by $\kappa$ for short. Assume that the following isomorphism
$\Psi: \nu(g) \cong k\kappa$ is well-defined, where $\nu(g)$ is the normal bundle of $g$, $k\kappa$ is the Whitney sum
of $k$ copies of the bundle $\kappa$. The triple $(g,\kappa,\Psi)$ represents an element in the cobordism group $Imm^{sf}(n-k,k)$, for a detailed definition and elementary properties see  \cite{A-F}. 

Assume that $k = \frac{n}{2}$, define the characteristic number $\theta(g,\kappa,\Psi) \in \Z/2$ as the parity
of self-intersection points of the immersion $g$ (we assume that $g$ is self-transversal). In fact,  $\theta$ depends not of $(\kappa, \Psi)$ and  $\theta(g,\kappa,\Psi) =\theta(g)$.

For given $k_-, k_+$, which satisfy the inequalities  $n> k_+ > \frac{n}{2}> k_->0$ the following sequence of homomorphisms
is well-defined:
\begin{eqnarray}\label{1}
Imm^{sf}(k_+,n-k_+) \stackrel{J^{k_+}_{\frac{n}{2}}}{\longrightarrow} Imm^{sf}(\frac{n}{2},\frac{n}{2}) \stackrel{J^{\frac{n}{2}}_{k_-}}{\longrightarrow}  Imm^{sf}(k_-, n-k_-).
\end{eqnarray}

\begin{definition}
Denote the subgroup $Im(J^{k_+}_{\frac{n}{2}}) \cap Ker(J^{\frac{n}{2}}_{k_-}) \subset Imm^{sf}(\frac{n}{2},\frac{n}{2})$ by $K^n_{k_+,k_-} \subset Imm^{sf}(\frac{n}{2},\frac{n}{2})$.
The 2-parameter family of subgroups determines the double filtration of the group $Imm^{sf}(\frac{n}{2},\frac{n}{2})$. (Note that $K^n_{\frac{n}{2},\frac{n}{2}}$ coincides with the whole group  $Imm^{sf}(\frac{n}{2},\frac{n}{2})$.)
\end{definition}

\begin{problem}{Generalized Kervaire Problem}\label{problem2}

For which $(k_+,k_-)$ the restriction of the homomorphism $\theta$ on the subgroup $K^n_{k_+,k_-}$ is non-trivial? 
(there exists $[(g,\kappa,\Psi)] \in K^n_{k_+,k_-} \subset  Imm^{sf}(\frac{n}{2},\frac{n}{2})$, for which
$\theta(g)=1$?)
\end{problem}  

\begin{proposition}\label{prop2}
A positive solution of the Kervaire Problem for 
$n=2^l-2$ is equivalent to a positive solution of Generalized Kervaire Problem for $K^n_{n-1,\frac{n}{2}}$.
\end{proposition}

\subsubsection*{Proof of Proposition $\ref{prop2}$}
Consider an immersion $f: M^{n-1} \looparrowright \R^n$.  Consider a submanifold  $M_1^{\frac{n}{2}} \subset M^{n-1}$, which is dual to the characteristic class $\kappa^{\frac{n}{2}-1} \in H^{\frac{n}{2}-1}(M^{n-1};\Z/2)$, in this case $\kappa: M^{n-1} \to \RP^{\infty}$ is the oriented class. (Below we use $\Z/2$-coefficients, in the case the coefficients are omitted.)  Consider the immersion  $g: M_1^{\frac{n}{2}} \looparrowright \R^n$, $g = f \vert_{M_1^{\frac{n}{2}}}$. Define the characteristic class  $\kappa \vert_{M_1^{\frac{n}{2}}}$, 
which is not re-denoted. By the construction a skew-framing  $\Psi: \nu_{g} \cong \frac{n}{2} \kappa$ is well-defined. 
The triple $(g,\kappa,\Psi)$ determines an element in the group $Imm^{sf}(\frac{n}{2},\frac{n}{2})$. 
The following equation is satisfied: $\Theta(f) = \theta(g)$. By the construction, $(g,\kappa,\Psi) \in K^{n}_{n-1,\frac{n}{2}}$. If the Kervaire Problem admits a positive solution, then there exists $f$, $\Theta(f)=1$. Then $\theta(g)=1$ and the Generalized Kervaire Problem for the term $K^n_{n-1,\frac{n}{2}}$ admits a positive solution. 

Conversely,  assume $\theta(g,\kappa,\Psi) = 1$. The inverse image of the subgroup $K^n_{n-1,\frac{n}{2}}$ in $Imm^{sf}(n-1,1)$by the homomorphism $J^{n}_{n-1,\frac{n}{2}}$  contains en element
$[f]$ with $\Theta([f]) = \theta([g,\kappa,\Psi]) = 1$. This proves that Kervaire Problem has the positive solution. \qed

\begin{proposition}\label{3}
Generalized Kervaire Problem for the term $K^n_{\frac{n}{2},\frac{n}{2}-1}$ is  positively solved, if and only if  $n=2,6,14$.
\end{proposition}

\subsection*{Proof of Proposition $\ref{3}$}

The group $K^n_{\frac{n}{2},\frac{n}{2}-1}$ is represented by triples $(g,\Xi)$, where
$g: M_1^{\frac{n}{2}} \looparrowright \R^n$, $\Xi$--is the trivialization of the normal bundle of $g$:
$\Xi: \nu_g \cong \frac{n}{2}\varepsilon$. Therefore, $(g,\Xi)$ represents an element in  $\Pi_{\frac{n}{2}}$. The value  $\theta(g,\Xi)$ coincides with the stable Hopf invariant $h([g,\Xi])$. By the Adams theorem $h([g,\Xi])=0$, 
if $l \ge 4$;  the homomorphism
$h: \Pi_{\frac{n}{2}} \to \Z/2$ is onto, if $n=2,6,14$. \qed

\section{Strong Kervaire Problem}

Define the element $w_n \in \pi_{2n-1}(S^n)$ by the Whitehead square of the generator 
$i_n \in \pi_n(S^n)$, $n \equiv 1\pmod{2}$. 

\begin{problem}{Strong Kervaire Problem}  
Assume $n=2^l-1$. For which $l$ the element $w_n$ is halved?
\end{problem}

\begin{theorem}\label{th1}
Strong Kervaire Problem is positively solved (for $n=2^{l+1}-2$ the element
$w_{\frac{n}{2}}$ is halved) iff Generalized Kervaire Problem for the term 
$K^{n}_{\frac{n}{2},\frac{n}{2}-2}$ is positively solved. 
\end{theorem}

\subsection*{Prof of Theorem $\ref{th1}$}
Denote $\frac{n}{2}$ by $k$ for short. By the	
Pontryagin Theorem the element in the homotopy group  $\pi_{2k-1}(S^k)$ is represented by a framed submanifold 
$(P^{k-1};\Xi_P)$, $P^{k-1} \subset S^{2k-1}$, $\Xi: \nu_P \cong k\varepsilon$. 
The element $w_k \in  \pi_{2k-1}(S^k)$ is represented by the Hopf link of the 2 standard framed spheres 
$(S^{k-1} \cup S^{k-1}) \subset \R^{2k-1} \subset S^{2k-1}$. Consider the standard internal embedding $I_1: S^{2k-1} \subset \R^n$, $n=2k$ and the concentric external embedding $I_2: S^{2k-1} \subset \R^n$.
By the assumption $w_k$ is even, therefore there exists a framed submanifold  $(Q^{k-1},\Xi_Q)$ in the internal sphere $I_2(S^{2k-1})$, which is the union of 2 disjoint copies of a framed embedding. Define an skew-framed immersion of a non-oriented manifold 
$g:N^{k} \subset \R^n$ with the only self-intersection point inside the sphere $I_1(S^{2k-1})$. 
The immersion  $g$ coincides with the Hopf link of the 2 standard framed spheres in the hypersphere $I_1(S^{2k-1})$, and coincides with $(Q^{k-1},\Xi_Q)$
in  the  hypersphere $I_1(S^{2k-1})$. In the exterior domain and in the middle  domain  $\R^n \setminus \{ I_1(S^{2k-1}) \cup I_2(S^{2k-1}) \}$
the immersion is a framed embedding.
In the interior ball the immersion is given by the pair of disks with the self-intersection point in the center of the ball.  In the middle domain the framed embedding is well-defined by the assumption that $w_n$ is halved,  in the exterior domain the framed embedding is given by the cylinder, which is bounded by the two copies of embedded submanifold with opposite framings (and with the opposite orientation).

By construction,  the oriented class  $\kappa_N: N^k \to \RP^{\infty}$ is an integer, $[\kappa_N] \in H^1(N^k;\Z)$. 
We get:
$\nu_g \cong k\kappa_N$, where the triple $(g, \kappa_N, \Psi_N)$ represents an element in  $Imm^{sf}(k,k)$, 
which is in the subrgoup $K^n_{\frac{n}{2},\frac{n}{2}-2}$, and $\theta(g)=1$.

The opposite is true, an arbitrary element  in  $K^n_{\frac{n}{2},\frac{n}{2}-2}$ is represented by a triple $(g, \kappa_N, \Psi_N)$, for which the characteristic class  $\kappa_N$ is integer. It is easy to check, that 
in a regular cobordism class there exists a triple $(g, \kappa_N, \Psi_N)$, where
$g$ is as described: $g(N^k)$ is divided by the spheres
$I_1(S^{2k-1})$, $I_2(S^{2k-1})$ into 3 submanifolds with boundaries.

Theorem  $\ref{th1}$ is proved.
\[  \]

\section{Barratt-Jones-Mahowald Theorem; James Theorem; Hill-Hopkins-Ravelel Theorem; Snaith Conjecture}

Let us formulate results, concerning Strong Kervaire Problem. 
The following theorem is proved in \cite{B-J-M}
\begin{theorem}\label{B-J-M}
Assume that for a given $n=2^{l+1}-2$ Strong Kervaire Problem is positively solved (this means that the homomorphism 
$\theta: K^{2^{l+1}-2}_{2^l-1,2^{l}-3} \to \Z/2$ is onto). Then Kervaire Problem for 
$n'=2^{l}-2$ is positively solved (this means that the homomorphism $\Theta_l: Imm^{sf}(2^l-3,1) \to \Z/2$ is onto). (Additionally, there exists an element of the order 2 in the stable homotopy group $\Pi_{2^l-2}$ of spheres, which is covered by Kahn-Priddy mapping  $\delta: Imm^{sf}(2^l-3,1) \to \ _2\Pi_{2^l-2}$ 
by an element $x \in Imm^{sf}(2^l-3,1)$ with $\Theta(x)=1$.)
\end{theorem}

Let us formulate the James Theorem \cite{Jam}.
Denote by $V_{k,2}$ the Stiefel manifold $V_{k,2}$ of $2$-framed in  $\R^k$;
the most interesting case is $k=2^l-1$, $ l \ge 2$. Denote by $I: V_{k,2} \to V_{k,2}$ the involution
by the formula: $\{e_1,e_2\} \mapsto \{e_1,-e_2\}$. 
\begin{definition}
Let us say that the manifold $V_{k,2}$ is neutral, if the involution $I$ is homotopic to the identity. 
\end{definition} 
  
\begin{theorem}[James Theorem]\label{Jm}
The manifold $V_{k,2}$, $k=2^l-1$, $ l \ge 2$ is neutral, iff the element  $w_k \in \pi_{2k-1}(S^k)$ is halved. 
This condition is equivalent to a positive solution of Strong Kervaire Problem for $n=2^l-2=k-1$.
\end{theorem}

\subsubsection*{Proof of Theorem $\ref{Jm}$}
Theorem 22.6 p. 139-141 \cite{J}.

\begin{example}\label{ex}
The manifold $V_{k,2}$ is neutral for $k=7$.
\end{example}

\subsection*{Proof of Example $\ref{ex}$}
Define a homotopy $F(t): V_{k,2} \to V_{k,2}$, $F(0)=Id$, $F(1)=I$. Consider the standard inclusion 
$\R^7 \subset \R^8$, which is orthogonal to the base vector $\1 \in \R^8$. Let $\f_1 \in \R^7 \subset \R^8$ be a base vector, which is orthogonal to the vector $1 \ne \f_1$. Denote by $L(\f_1)$ the orthogonal complement to the vector  $\f_1$ in the subspace $\R^{7}$. By construction, $\e_2 \in (\e_1)^{\perp}$.

A one-parameter family of orthogonal transformation 
$G(t,\e_1): \R^7 \to \R^7$  which transforms the subspace $(\e_1)^{\perp}$ to itself, is well-defined by multiplication of vectors from $(\e_1)^{\perp}$ by the  unite  ${\bf 1} \cos(t\pi) + \e_1 \sin(t\pi)$. Denote by $F(t,\e_1): \R^7 \to \R^7$ 
the transformation, which is the identity in the subspace generated by the vector  $\e_1$ and which coincides to $G(t,\e_1)$ on the subspace $(\e_1)^{\perp}$. The homotopy $F(t,\e_1)$ is required.
\[  \]
\subsubsection*{Remark}
Example $\ref{ex}$ could be interesting with respect to results \cite{C-G-K-Sh}. The space $V_{k,2}$ is the model
of random closed $k$-gones on the plane with the total unite length of edges. The homotopy $F(t)$ determines 
the transformation of the space of $k$-gones to itself,  a gone is transformed to its mirror image. 
\[  \]

It is known, that Kervaire Problem  is positively solved for
$n=2,6,14,30,62$.  In the case $n \ge 254$ Kervaire Problem is negatively solved by Hill-Hopkins-Ravenel Theorem, \cite{H-H-R}. A negative solution of Kervaire Problem except, probably, a finite values of $n$, using different arguments,  was presented in 2011 Workshop on Kervaire invariant and stable homotopy theory, Edinburgh, by the author. 
(Let me remark that Lemma 37 in \cite{A} is not well proved; Lemma 36 is sufficient.)

The following conjecture is open \cite{S}.

\begin{conjecture}{Snaith Conjecture}\label{Sn}
Kervaire Problem is negatively solved for 
 $n=126$.
\end{conjecture}
\[  \]

If Snaith Conjecture is not true, Strong Kervaire Problem for
 $n+1=127$ is solved negatively.

\section{Solution of Strong Kervaire Problem for $\frac{n}{2}=15$}

Define a closed manifold 
 $M^{15}$ as a semi-direct product of $S^7 \times S^7$ and the circle $S^1$, which is fibered over $S^1$ with the
 fiber $S^7 \times S^7$. The identification of the fiber $S^7 \times S^7$ along $S^1$ is given by the 
 involution 
 $S^7 \times S^7 \to S^7 \times S^7$, which permutes the factors. The orientation class $w_1(M^{15})=\kappa \in H^1(M^{15};\Z/2)$ is induced from the generator $H^1(S^1)$ by the projection $M^{15} \to S^1$.

Let us define an immersion
$\varphi: M^{15} \looparrowright \R^{30}$, which is skew-framed  by $\Psi: \nu_{\varphi} \cong 15\kappa$.  Consider
the manifold  $\hat{K}^{15}=\RP^7 \times \RP^7 \tilde{\times} S^1$, which is fibered over $S^1$ with the fiber $\RP^7 \times \RP^7$.  The identification of the fiber $\RP^7 \times \RP^7$ along $S^1$ is given by the 
 involution 
 $\RP^7 \times \RP^7 \to \RP^7 \times \RP^7$, which permutes the factors.  
 
The normal bundle $\nu_{\hat{K}}$ of the manifold 
 $\hat{K}^{15}$ is isomorphic to the Whitney sum of an odd number $s$ of copies of the line orienting bundle over
$\hat{K}^{15}$, denote this orienting bundle by  $\hat{\kappa}$. We may assume that $s=15$, because 
$2\hat \kappa \cong 2\varepsilon$ over $\nu_{\hat{K}}$, the normal bundle can be replaced by 
$\nu_{\hat{K}} \cong (15-2r)\hat \kappa \oplus 2r\varepsilon$.

An immersion $\hat \varphi: \hat K^{15} \looparrowright \R^{30}$ is well-defined, this immersion is skew-framed by
$\hat \Psi: \nu_{\hat \varphi} \cong 15 \hat \kappa$. Consider the standard 4-sheeted covering 
$p: M^{15}=(S^7 \times S^7) \tilde{\times} S^1 \to \hat{K}^{15}$, which is induced by the pair of coverings $S^7 \to \RP^7$. The immersion  $\varphi$ is well-defined by the formula:  $\varphi = \hat \varphi \circ p$. A skew-framing $\Psi$ is well-defined by the formula: $\Psi = \varphi^{\ast}(\hat \Psi)$.
The value $\theta(\varphi,\kappa,\Psi)$ is calculated as a parity of self-intersection points of the immersion $\varphi=\hat \varphi \circ p$, which is generically deformed into  $\R^{30}$ by a small perturbation. Let us prove the formula $\theta(\varphi,\kappa,\Psi)=1$.

To calculate 
$\theta \in \Z/2$, we describe a self-intersection manifold of the (non-generic) immersion $\hat \varphi \circ p: M^{15} \looparrowright \R^{30}$, the self-intersection manifold contains $2$ components.  A component $N_1^{15}$ consists of points $[(x_1,x_2;y_1,y_2]$, where $x_1=x_2$, 
$y_1 = Ty_2$, or $Tx_1=x_2$, 
$y_1 = y_2$, where $T: S^7 \to S^7$  is the standard involution on the factors, $x_1, x_2$ are points on the same factor
in the (local) coordinates system $M^{15} \cong S^7 \times S^7 \times z$, $z \in S^1$. The component  $N_1^{15}$ is a semi-direct product of a non-connected manifold  $L_1^{14}$ with the circle. 
Obviously, after a regular generic perturbation of $\hat \varphi \circ p$, the component gives no 
contribution to the value $\theta$.

The last component denote by 
 $N^{15}$, this component consists of points $[(x_1,x_2;y_1,y_2]$, where $x_1=Tx_2$, 
$y_1 = Ty_2$. The manifold is also a semi-direct product of a 14-manifold, denoted by $L^{14}$, with the circle $S^1$. 
The characteristic class of the canonical covering $\tilde N^{15} \to N^{15}$ over self-intersection manifold denote by
 $p_N \in H_1(N^{15})$. The normal bundle
$\nu_N$ over $N^{15}$ in $\R^{30}$ is isomorphic to the Whitney sum $\kappa_N \oplus 14\varepsilon$, where $\kappa_N$
is the linear bundle with the generator along the circle $S^1$. To calculate $\Theta$  the following 
equation is required:
\begin{eqnarray}\label{44}
\langle p_N^{14} (p_N + \kappa_N);[N] \rangle =1
 \end{eqnarray}

Describe explicitly the manifold $N^{15}$ and the submanifold $L^{14} \subset N^{15}$;
describe the projection $\pi_N : N^{15} \to (\RP^7 \times \RP^7) \times S^1$
(by the restriction  $\pi_N \vert_L$ of the projection the submanifold $L^{14}$ is mapped onto $(\RP^7 \times \RP^7) \subset (\RP^7 \times \RP^7) \tilde \times S^1$).

The manifold $L^{14}$ are defined by two  non-ordered pairs of points 
$[(x_1; x_2; y_1; y_2]$,   $S^7 \times S^7$, each point in the pair $(y_1,y_2)$ belongs to the second factor, and $T(x_1) = x_2, T(y_1) = y_2$.  Describe the canonical covering $p_N :
\tilde N^{15}\to N^{15}$ in terms of submanifold $(\RP^7 \times \RP^7) \tilde \times S^1$. 

Denote by $\pi_L$ the restriction of $p_N$ over the submanifold $L^{14}$. The class 
$\pi_L = t_1 + t_2 \in H_1(\RP^7 \times \RP^7)$ is the sum of the generators of the two factors. 

Alternatively, the same characteristic class is defined by the formula: 
$\pi_L = f^{\ast}(t)$, $t \in H^1(\RP^7)$, where
$f : \RP^7 \times \RP^7 \to \RP^7$ is the mapping, determined by the formula:  $f(a; b) =
a \odot b$, where $a, b$ are points on the corresponding factor $\RP^7$, $\odot$ is the product in the Cally
algebra. 

The manifold $L^{14}$ coincides with the  2-sheeted covering over  $\RP^7 \times \RP^7$ with the characteristic class
$\pi_L \in H^1(\RP^7 \times \RP^7)$. Obviously, 
the class $\pi_L$ is the restriction of a class  $\pi_N \in
H^1(\RP^7 \times \RP^7) \tilde \times S^1)$ by the inclusion $(\RP^7 \times \RP^7) \subset (\RP^7 \times \RP^7)
\tilde \times S^1$.

The class $\pi_N$ is determined by the formula:
$\pi_N = g^{\ast}(t)$, $g : (\RP^7 \times \RP^7) \tilde \times S^1 \to \RP^7$,
$g(a; b; z) = a \odot b$, where a triple $(a; b; z)$ represents a point on $(\RP^7 \times \RP^7) \tilde \times S^1$.
The manifold $N^15$ coincides with the covering space of the 2-sheeted covering over $(\RP^7 \times \RP^7) \tilde \times S^1$ with characteristic class $\pi_N$.

Describe the characteristic class
$p_N \in H^1(N^{15})$ of the canonical covering  $\tilde N^{15} \to
N^{15}$. Firstly, describe the class $p_L \in H^1(L)$, which is defined as the restriction of $p_N$ 
on the submanifold $L^{14} \subset  N^{15}$. Consider the projection 
$s_1 : L^{14} \to \RP^7$, $s_1 = p_1 \circ \pi_L$, $p_1 : \RP^7 \times \RP^7 \to \RP^7$ is the projection
on the first factor. Consider the projection  $s = p_2 \circ \pi_L$, where
$p_2 : \RP^7 \times \RP^7 \to \RP^7$ is the projection on the second factor. The mappings
 $I \circ s_1 : L^{14} \to \RP^7 \subset \RP^{15}$, 
$I \circ s_2 : L^{14} \to \RP^7 \subset \RP^{15}$ are homotopic, where $I : \RP^7 \subset \RP^{15}$
is the standard inclusion. In particular, $s^{\ast}_1(t) = s^{\ast}_2(t)$. Therefore, the mapping  $p : N^{15} \to
\RP^{15}$ is well-defined (this mapping represents the class $p_N \in H^1(N^{15})$) by a gluing of the two mappings
$s_{+\varepsilon} = I \circ s_1 : L^{14} \times \{+\varepsilon\} \to \RP^{15}$, 
$s_{-\varepsilon} = I \circ s_2 : L^{14} \times \{ -\varepsilon\} \to \RP^{15}$
(this two mappings give boundary conditions on the preimage of a short segment  $S^1 \setminus (-\varepsilon,+\varepsilon) \subset S^1$
by the projection $\kappa_N : N^{15} \to S^1$); after the gluing we get the mapping of a short segment 
$L^{14} \times [-\varepsilon,+\varepsilon] \subset N^{15}$ with the central submanifold  $L^{14} \cong L^{14} \times \{0\} \subset
N^{15}$, which we denote by  $s_[-\varepsilon,+\varepsilon] : L^{14} \times  [-\varepsilon,+\varepsilon] \to \RP^{15}$. The following formula is required:
\begin{eqnarray}\label{33}
\langle p_N^{15};[N] \rangle =1.
\end{eqnarray}
This formula implies that the the mapping $p : N^{15} \to \RP^{15}$, which represents the  class $p_N$, 
satisfies the condition: $deg(p) = 1 \pmod{2}$.

The formula $(\ref{33})$ is deduced from Lemma $\ref{lem}$ below. Consider the mapping 
$s_[-\varepsilon,+\varepsilon]$ with prescribed boundary conditions  $s_[-\varepsilon,+\varepsilon] : (L^{14} \times
[-\varepsilon,+\varepsilon]; \partial) \to (\RP^{15},\RP^7)$. Let us prove that the degree
of this mapping (the number of regular preimages) is odd.

\begin{lemma}\label{lem}
The following formula is satisfied: $deg(s_[-\varepsilon,+\varepsilon]) = 1$.
\end{lemma}

\subsection*{Proof of Lemma $\ref{lem}$}
In the standard projective $15$-space  $\RP^{15}$ consider the two standard disjoint subspaces
 $\RP^7_- \cup \RP^7_+ \subset \RP^{15}$, this two projective subspaces are diffeomorphic: $\RP^7_- \cong \RP^7_+$.
Denote by $U_{-} \cup U_{+} \subset  \RP^{15}$ its thin regular neighborhoods. Note that
the manifold $L^{14} \times [-\varepsilon,+\varepsilon]$ is diffeomorphic to $\RP^{15} \setminus (U_{-} \cup U_{+})$.
The diffeomorphism $L^{14} \times [-\varepsilon,+\varepsilon] \cong \RP^{15}  \setminus (U_{-} \cup U_{+})$
is satisfied the following condition: the two projections 
$pr_{-} : \partial U_{-} \to \RP^7_-$,
$pr_{+} : \partial U_{+} \to \RP^7_+$ of the boundary components 
onto the central submanifold coincides with the mappings 
$s_{-\varepsilon} : L^{14} \times \{-\varepsilon\} \to  \RP^{7}_-$,
$s_{+\varepsilon} : L^{14} \times \{+\varepsilon\} \to  \RP^{7}_+$
by this diffeomorphism. To conclude the proof, the homotopy between the mappings  $pr_-$, $pr_+$ is considered, 
the image if this homotopy 
is inside $8$-dimensional subspace, this subspace contains no the center of the bottom cell in 
 $\RP^{15}$. Lemma  $\ref{lem}$ is proved.

Calculate the class $(\ref{33})$. More precisely, we prove that the equations  $(\ref{33})$, $(\ref{44})$
are equivalent. This follows from the equation: 
$\langle p_N^{14} \kappa_N;[N] \rangle =0$. The last equation is obvious, because 
 $\langle p_N^{14} \kappa_N;[N] \rangle = \langle p_L^{14};[L] \rangle = 0$.
The class $p_L$ is represented by a mapping  $L^{14} \to \RP^{\infty}$ with the image inside 
the $7$-skeleton (therefore the class is a coboundary). It is proved that a small regular alteration of the mapping  $\hat \varphi \circ p$  self-intersects by an odd number of points.

\section{ The Jones Manifold}

In the paper \cite{J} a framed closed manifold $\tilde M^{30}$ of dimension  $30$ with Arf-invariant $1$
is constructed. 

\subsubsection*{Definition of the quadratic form $q: H_{2k+1}(\tilde M^{4k+2}) \to \Z/2$ for a framed $4k+2$-manifold}

Let $(\tilde M^{4k+2},\Xi)$ be an arbitrary closed manifold with a prescribed trivialization of the stable normal bundle 
$\Xi: \nu_{\tilde M} \cong r \varepsilon$, $r-2 \ge 4k+2 = dim(\tilde M)$. The framing  $\Xi$ determines up to regular concordance an immersion $\varphi: (\tilde M^{4k+2}\setminus pt) \looparrowright \R^{4k+2}$. 
Let $x \in H_{2k+1}(\tilde M^{4k+2})$ be a homology class. By Thom Theorem and Whitney Theorem the class $x$
is represented by an embedding
 $l: K^{2k+1} \subset \tilde M^{4k+2}$, where $K^{2k+1}$ is a closed, generally speaking, non-oriented closed $2k+1$-manifold. 

Consider the immersion
$\varphi \circ l: K^{2k+1} \subset (\tilde M^{4k+2}\setminus pt) \looparrowright \R^{4k+2}$. Define $q(x,l)$ as
the parity of self-intersection points of the immersion  $\varphi \circ l$.
The value $q: H_{2k+1}(\tilde M^{4k+2};\Z/2) \to \Z/2$ depends only on $x$ and depends not on 
the cycle
 $l:K^{2k+1}\subset \tilde M^{4k+2}$, which represents $x$. The form $q$ is quadratic and not degenerated, because this form is associated with the intersection of $2k+1$-cycles on $\tilde M^{4k+2}$. Define the Arf-invariant of the form $q$ by $Arf(q) \in \Z/2$. It is not hard to how, that 
$Arf(q)$ is well-defined and also is invariant with respect to framed surgery of the framed manifold  $(\tilde M^{4k+2},\Xi)$. By Pontryagin theorem $Arf(q)$ determines a homomorphism $\Pi_{30} \to \Z/2$. 

Assume that the manifold
 $\tilde M^{4k+2}$ is  $2k-1$-connected and $k \cong 3 \pmod{4}$, $4k+2 \ge 22$. In this case the definition of $q: H_{2k+1}(\tilde M^{4k+2}) \to \Z/2$ coincides with a well-known definition from  \cite{Br}. Namely, assume that the homology class  $x$ is represented by an embedded sphere $l: S^{2k+1} \subset M^{4k+2}$, one define $q(x)=1$, if the normal bundle $\nu(l)$ of the embedding $l$ is non-trivial, and define  $q(x)=0$, if the normal bundle of $l$ is trivial. Let us deduce this definition from the definition of $q$, given above. The normal bundle $\nu(l)$
is isomorphic to the normal bundle of the immersion
$\varphi \circ l: S^{2k+1} \looparrowright \R^{4k+2}$, the number of self-intersection points of this immersion 
is odd iff  the bundle $\nu(l)$ is non-trivial. 
\[  \]


Assume $4k+2=30$. Let us construct a stably parallelized manifold $\tilde M^{30}$, which is called the Jones Manifold  \cite{J}. Denote by  $P^2$  the closed non-orientable surface of the Euler characteristic  $+3$
(the connected sum of the projective plane and the torus). Consider an arbitrary embedding $i_P: P^2 \subset \R^{35}$, denote by $\nu_P$ the normal bundle of the embedding $i_P$, $\dim(\nu_P)=33$. Denote by  $\nu_0 \subset \nu_P$ a trivial subbundle of the dimension $32$, $\nu_0 \cong 32\varepsilon$. Let us write $\nu_P \cong \hat \nu_P \oplus \nu_0$, where $\hat \nu_P$ is the line orienting normal bundle over  $P$. Decompose the fiber of the bundle
$\nu_0$ 
over the marked point $pt \in P^2$ into  4 blocks, denote the blocks by $A,B,C,D$, the each 
block is 8-dimensional.

Consider the decomposition of the surface $P^2 \cong Q^2 \cup_S R^2$, where $Q^2$ is the projective plane
with hole, $R^2$ is the torus with hole, the two surfaces $Q^2$, $R^2$ are glued along the common boundary circle $S$. 
Denote by
$a \in H_1(Q^2)$ the generator on the projective plane, by $b_1, b_2 \in H_1(P^2)$  generators on the torus.
Define a flat 4-bundle $\mu: E(\mu) \to  P$, the structure group of this bundle is given by
the permutation of the $1$ and $3$ base vectors along the loop on $Q^2$, represented $a$;
by the permutation $(1,2)(3,4)$ of the base vectors along the loop, represented  $b_1$; and by the permutation
$(2,3)(4,1)$ of the base vectors along the loop, represented $b_2$. The bundle $\mu$ has the dihedral structure group $\D$, this group contains $8$ elements, the formulas determine a representation 
$\mu: \pi_1(P^2) \to O(4)$.

Consider the bundle $8\mu$, the Whitney sum of $8$ copies of the bundle  $\mu$. The bundle $8\mu$ over $P^2$ is trivial,  therefore an isomorphism  $8\mu \cong \nu_0$ is well-defined, fibers of the bundle $8\mu$ (in, particular, the fiber over the marked point $pt \in P$) are decomposed into  $8$-bundles as follows: $A \oplus B \oplus C \oplus D$.
Obviously, $\mu \cong \varepsilon \oplus \hat\mu$, where $\dim(\hat \mu)=3$. It is not hard to prove, that $w_1(\hat \mu) \in H^1(P^2)$ is dual to  $a$, $w_2(\hat \mu) \in H^2(P^2)$ represents the fundamental class on $P^2$.
Therefore we have
$\hat \Xi: \hat \nu_P \oplus \hat \mu \cong 4\varepsilon$ and
$\Xi: \hat \nu_P \oplus \mu \cong 5 \varepsilon$.

Let us consider the fiber  $A \oplus B \oplus C \oplus D$ of the bundle $\nu_0$ over the marked point $pt \in P^2$ and
consider the Hopf immersion $h: S^7 \looparrowright \R^8$, which determines a generator of $[\Pi_7]_{(2)} \cong \Z/16$.
Consider the Cartesian product $h^4: \looparrowright A \oplus B \oplus C \oplus D$ of $4$ copies of the immersion $h$.
Define the required manifold $\tilde M^{30}$ as the semi-direct product  $P^2 \tilde \times (S^7)^4$, where
the dihedral group $\D$ transforms the fiber  $(S^7)^4$ by the formulas of $\mu$ for generators $a, b_1, b_2$. By
the construction, $\D$ keeps the immersion $h^4$ and keeps the framing if this immersion. 

The immersion 
 $\tilde M^{30} \looparrowright E(8\mu)$ into the total space of the bundle  $8\mu$ is well-defined. Because $E(8\mu)$ is diffeomorphic to the total space of the bundle $\nu_0$, the immersion 
$f: \tilde M^{30} \looparrowright \R^{35}$ is well-defined. Obviously, the normal bundle  $\nu_f$ of the immersion $f$ is isomorphic to the Whitney sum $\hat \nu_P \oplus \mu$; by the isomorphism $\Xi$ we get the trivial normal bundle $\nu_f$. The manifold $\tilde M^{30}$ is framed and  $f$ is a framed immersion, which determines an element in   $\Pi_{30}$. Jones Manifold is well-defined. 

\subsection*{Calculation of Hamiltonian base in the middle homology of the Jones Manifold}
 
Let us prove that  the vector space $H_{15}(\tilde M^{30})$ is  $8$-dimensional, then determines the Hamiltonian base of this space with respect to the standard bilinear product. 
Take the marked point  $pt \in P^2$ and marked the $7$-cycles in the fiber $(S^7)^4$ of the projection $p: \tilde M^{30} \to P^2$ by $A,B,C,D$ correspondingly with the subfibers of the vector bundle  $\nu_0$.   
Cycles $A \otimes B$, $A \otimes C$, $A \otimes D$, $B \otimes C$, $B \otimes D$, $C \otimes D$ are the base cycles, in $H_{14}((S^7)^4)$, $(S^7)^4 = p^{-1}(pt)$. Denote $H_{14}((S^7)^4)$ by $\Omega$ for short.

The punctured surface $P^2 \setminus pt$ is homotopy equivalent to a bouquet of 3 circles:  
 $P^2 \setminus pt \cong S^1_a \vee S^1_{b_1} \vee S^1_{b_2}$, where the circles  $S^1_a, S^1_{b_1}, S^1_{b_2}$
are marked correspondingly to cycles in $H_1(P^2)$.  Along the each circle a monodromy is well-defined:
$\theta_a: \Omega \to \Omega$, $\theta_{b_1}: \Omega \to \Omega$, $\theta_{b_2}: \Omega \to \Omega$.
The each monodromy admits $2$-dimensional fixed subspace, to the complement of this subspace
the monodromy is represented by a pair of permutations of base vectors.

We may take the denotations of the cycles, such that
$[A \times B \times b_1], [C \times D \times b_1]$ are closed $15$-cycles over $b_1$;  
$[B \times C \times b_2], [D \times A \times b_2]$ are closed  $15$-cycles over $b_2$.  
This means that the monodromy along $b_1$ transposes the 2 pairs of the cycles $A$ and $B$, $C$ and $D$.
Therefore along $b_1$ there exists the extra two invariant cycles  $A \times C + B \times D$,
$A \times D + B \times C$. The 2 pairs of cycles  ($A \times C$ and $B \times D$), ($A \times D$ and $B \times C$) are transposed along $b_1$. Analogously, the monodromy along $b_2$ is described. 

Assume that $7$-cycles $A,C$ are invariant with respect to the monodromy $\theta_a$ along $a$.
Therefore, the cycles  $A \times C$, $B \times D$,
$A \times B + A \times D$, $C \times B + C \times D$ are invariant. The 2 pairs of cycles ($A \times B$ and $A \times C$), $(C \times B)$ and
$C \times D$) are transposed by $\theta_a$.
 
This gives the dimension of the $\Z/2$--vector space: $\dim H_{15}(\tilde M^{30} \setminus p^{-1}(U(pt))) = 12$. 
Then let us calculate the dimension of the image of the homomorphism $\partial_{\ast}: H_{15} (\partial U(pt)) \to H_{15}(\tilde M^{30} \setminus $$p^{-1}(U(pt)))$,  which induced by the inclusion of the boundary into the manifold with boundary. It is not hard to prove that  the image of this homomorphism is $4$-dimensional.
The dimension 
$\dim(H_{15} (\partial U(pt)))=6$, the $2$ base cycles are in the kernel  $Ker(\partial_{\ast})$:
namely, the cycles $[A \otimes B + C \otimes D] \otimes [\partial U]$, $[ B \otimes C + C \otimes D] \otimes [\partial U]$.
Therefore, we get $\dim(H_{15}(\tilde M^{30}))=8$.

Let us calculate the Hamiltonian base in $H_{15}(\tilde M^{30})$, this base contains 4 pairs of cycles, 
cycles in each pair are dual, cycles from different pairs are orthogonal.  By construction, the projection
$p: \tilde M^{30} \to P^2$ on the base with the fiber $(S^7)^4$ is well-defined. Denote by $\hat Q \subset \tilde M$, $\hat R \subset \tilde M$, $\hat S \subset \tilde M$ preimages by $p$ of corresponding surfaces and its common boundary. 

Assume that the marked point $pt \in S \subset Q^2$. The cycles 
 $[A \times C \times a], [B \times D \times a]$ determines a Hamiltonian pair. Consider the cycle  
$[A \times B \times b_1]$ and define a cycle  $[C \times D \times b_2]'$ as the corresponding cycle in the
preimage of the curve $ab_2a^{-1}$ (this curve is embedded into  $P^2$), which coincides with the product 
$C \times D$ over $pt \in P^2$.
Cycles $[A \times B \times b_1]$, $[C \times D \times b_2]'$ determines a Hamiltonian pair.  

Analogously, a Hamiltonian pair  
$[A \times B \times b_2]'$, $[C \times D \times b_1]'$ is well-defined.
The cycle $[A \times B \times b_2]'$ is defined as a cycle in the preimage of the curve
$a b_2 a^{-1}$, which coincides with the product $A \times B$ over $pt$.

Finally, the last Hamiltonian pair
$[A \times C \times b_1]'$, $[B \times C \times b_2]$ is analogously defined. In construction of the cycle $[B \times C \times b_2]$
we used that
the cycle $B \otimes C$ is invariant along the curve $b_2$, the cycle
$[A \times C \times b_1]'$ is the result of an extension of $A \otimes D$
along the curve $a b_1 a^{-1}$ from $pt$.

Let us prove that the coefficient of intersection of the cycles
$[A \times D \times b_1]'$, $[C \times D \times b_2]'$ is trivial (the last formula for intersection coefficients
are analogous).
The curves $ab_1a^{-1}$, $ab_2a^{-1}$ on $P^2$ intersects into $5$ points
(denote the intersection set by $X$). One of the point from $X$ is the intersection point of the curves $b_1$ and $b_2$, the last $4$ points are intersection points of a regular neighborhood of the cycle with its parallel copy.
By a direct consideration, in each intersection point self-intersection of cycles is canceled by a small deformation
in the fibers over the intersection points.

The Hamiltonian base in 
$H_{15}(\tilde M^{30})$ is well-defined.

\subsubsection*{Calculation of Arf-invariant of Jones Manifold}
Let us calculate the quadratic form $q$ for vectors in the Hamiltonian base of
$H_{15}(\tilde M^{30})$. Prove that $q([A \times C \times a])=0$; $q([A \times C \times b_1]')=1$, $q([B \times D \times b_2]')=1$. Values of $q$ on the last base vectors are correlated and give no contribution to Arf-invariant. 

The equation $q([A \times C \times a])=0$ is satisfied, because the cycle is represented by a framed 
$15$-manifold, $q([A \times C \times a])$ coincides with the Hopf invariant, by the Toda Theorem the Hopf invariant in dimension $15$ is trivial. 
The equation $q([A \times C \times b_1]')=1$ is proved using results of Section 5. The considered cycle
is represented by  the manifold
$(S^7 \times S^7) \tilde \times S^1$, this manifold is diffeomorphic to the manifold $\tilde M^{15}$, this manifold if 
immerser into $\R^{30}$ by
$\varphi: \tilde M^{15} \looparrowright \R^{30}$. By construction,
the restriction of the immersion $\varphi$ on the submanifold $S^7 \times S^7 \subset (S^7 \times S^7) \tilde \times S^1$ coincides with the immersion of the submanifold $\tilde L^{14}$, as in the example.

Using direct elementary consideration one may prove that the parity of self-intersection points of $\varphi$
is well-defined from this data.
It is proved that the framed manifold  $(\tilde M^{30}, \Xi)$ is of the $Arf$-invariant 1.

\section{Browder Theorem, Eccles Theorem, the positive solution of Kervaire Problem in dimension $n=30$}

An element
$\alpha$ in the stable homotopy group of spheres  $\Pi_{30} = \pi_{n+30}(S^{n}$, $n \ge 32$, is represented by a framed manifold 
$(\tilde M^{30},\Xi)$.  By the Browder theorem \cite{Br} a framed manifold  $(\tilde M^{30},\Xi)$ has Arf-invariant one 
iff the element $\alpha$ in $\Pi_{30}$ is detected by a secondary cohomological operation, which is based on the Adem relation
in the Steenrod algebda:

\begin{eqnarray}\label{adem1}
 Sq^{2^j}Sq^{2^j} + \sum_{i=1}^{j-1} Sq^{2^{j+1}-2^i} Sq^{2^i} = 0, \quad j=4.
\end{eqnarray}

By Eccles Theorem
 \cite{Ecc} the element $\alpha$ is detected by the considered cohomology operation, iff the element
 $\beta = \lambda_{\ast}(\alpha) \in \pi_{29+n}(\Sigma^n(K(\Z/2,1)))$, which is associated with  $\alpha$ by the Khan-Priddy mapping $\lambda: Q\RP^{\infty} \to Q(S^0)$,
is detected by the operation $Sq^{16}$ in the cone of the mapping  $\beta$ using the generic cohomology class in  $H_{n+15}(\Sigma^{n}(\RP^{\infty}))$. By Pontyagin-Thom-Wells Construction the element  $\beta$ is represented by
an immersion  $f: M^{29} \looparrowright \R^{30}$ in codimension 1 of a closed, generally speaking, non-orientable manifold. Moreover, the characteristic number $\Theta(f)=1$, iff  $Sq^{16}$ detects$\beta$. In the presented form the statement is used in \cite{A-E}. 

Because the Jones Manifold is a framed manifold of the Arf-invariant one, there exists an immersion
$F: M^{39} \looparrowright \R^{30}$ with $\Theta(f)=1$. Therefore Kervaire Problem is positively solved in dimension $n=30$. 

\section{Eccles-Wood filtration. Estimation of the filtration for a standard Kervaire $30$-manifold}

The fundamental paper  \cite{K} is dedicated to applications of homotopy theory in differential topology. 
In this paper M.Kervaire constructed a closed $4$-connected $PL$-manifold of dimension $10$, which admits no smooth structure.
In the dimension $30$ the analogical  $14$-connected manifold is smooth, called a standard Kervaire manifold.
A standard Kervaire manifold is defined as a result of a normal surgery of the Jones manifold
$\tilde M^{30}$ into a $14$-connected manifold $N^{30}$, 
$dim(H_{15}(\tilde N^{30};\Z/2))=2$. A standard Kervaire manifold is well-defined up to a connected sum with 
a homotopy $30$-sphere. By construction,  
$N^{30}$ is framed. It is well-known, that $\tilde N^{30}$ admits a $PL$-embedding into $\R^{32}$. 

\begin{theorem}\label{trivial}
There are no smooth embedding 
 $N^{30} \times D^{14} \subset \R^{46}$, where $N^{30}$ is a standard Kervaire manifold, $D^{14}$ is the standard disk.
\end{theorem}

To prove the theorem one has to estimate the Eccles-Wood filtration \cite{Ecc1} for a standard Kervaire manifold.
Let $N^n$ be a closed framed manifold (let us restrict to the case $n \equiv 0 \pmod{2}$), with a prescribed trivialization of the stable normal bundle: $\Xi: \nu_N \cong r\varepsilon$ of a dimension $r$, $r\ge n+2$. 
One say that 
$(N^n,\Psi)$ admits the Eccles-Wood filtration, at least   $(k_+,k_-)$, $n > k_+ > \frac{n}{2} \ge k_- > 0$, if
there exists an embedding  $N^n \subset \R^{2n - k_-}$ with the normal vector field $\psi$ of the dimension 
$k_+ - k_-$ (of the codimension  $n-k_+$). 

The following statement is equivalent to the theorem 
 $\ref{trivial}$.
\begin{proposition}\label{trivia}
The standard Kervaire manifold $(N^{30},\Xi)$ admits not the Eccles-Wood filtration 
$(28,14)$. 

\subsubsection*{Remark}
In \cite{A-C-R} a more stronger result is proved: the standard Kervaire manifold $(N^{30},\Xi)$ admits not the Eccles-Wood filtration 
$(23,14)$. 
\[  \]

\end{proposition}

\subsection*{Proof of Proposition $\ref{trivial}$}
Assume, there exists an embedding
$f: N^{30} \subset \R^{46}$, for which the normal bundle $\nu_f$ admits
$14$ linear independent sections. Prove that the immersion $f$ cannot be an embedding.
By the assumption, an isomorphism of the normal bundle 
 $\nu_f \cong \eta \oplus 14\varepsilon$, where $\eta$ is a real $2$--bundle, is well-defined. 
Because  
 $\pi_1(N^{30}) = \pi_2(N^{30})=0$, the bundle $\eta$ is a trivial bundle.
Therefore a framed immersion $(f,\Xi)$, $\Xi: \nu_f \cong 16\varepsilon$ is well-defined. 
It is well-known, that the Arf-invariant for a framed manifold
 $(N^{30},\Xi) \in \Pi_{30}$ depends not on a framing  \cite{J-R}. Therefore for the framed immersion 
$(f,\Xi)$ we get $Arf(f,\Xi)=1$.  
 
Let us consider a hyperprojection 
$\pi: \R^{46} \to \R^{45}$ and consider 
the mapping $g=\pi \circ f$, by the compression theorem \cite{R-S} this mapping could be defined as a framed immersion (by a framing $\hat \Xi$, where $\hat \Xi \oplus Id_{\varepsilon} \cong \Xi$), which represents an element $[(g,\Psi)])$, in the cobordism group $Imm^{sf}(15,15) \cong \Pi_{15}$. 

Consider the self-intersection manifold of the framed immersion 
$g: N^{30} \looparrowright \R^{45}$. By construction,  $[(g,\Psi)] \in K^{30}_{29,14}$
(a priori we have $[(g,\Psi)] \in K^{30}_{28,14}$). The first index $k_1$ of the filtration is the maximal,
therefore by Browder and Eccles Theorems the equation 
$\Theta([(g,\Psi)])=1$ is satisfied, because $Arf(g,\hat\Xi)=1$. The second index of the filtration
shows that 
$\Psi$ is a frame, by Proposition $\ref{3}$ we get $\Theta([(g,\Psi)])=0$. Therefore an embedding $N^{30} \times D^{14} \subset \R^{46}$ cannot exists. \qed
 
The author is grateful to Popelenskii Th. Yu. for discussions.

\end{document}